\documentclass[12pt]{amsart}
\usepackage{color}
\usepackage{ulem,ifthen,xcolor,xkeyval,pdfcolmk}

\usepackage[abbrev]{amsrefs}

            \definecolor{pink}{rgb}{1,0,1}

\renewcommand{\bar}{\overline}
\newcommand{\eps}{\epsilon}
\newcommand{\pa}{\partial}

\newfont{\fnt}{cmr10 scaled 550}
\renewcommand{\eps}{\varepsilon} 
\newcommand{\n}{\nabla}

\newtheorem{theorem}{Theorem}[section]

\newtheorem{conj}{Conjecture}
\newtheorem{lemma}{Lemma}[section]

\newtheorem{prop}{Proposition}[section] 

\newtheorem{definition}{Definition}[section]

\theoremstyle{remark}
\newtheorem{rem}{Remark}[section] 
\numberwithin{equation}{section} 
\newtheorem{question}{Question}[section]

\font\strange=msbm10

\renewcommand{\epsilon}{\varepsilon}

\renewcommand{\Sigma}{\varSigma}
\newcommand{\R}{{{\mathchoice  {\hbox{$\textstyle{\text{\strange R}}$}}
{\hbox{$\textstyle{\text{\strange R}}$}}
{\hbox{$\scriptstyle  N\kern-0.3em  R$}}  
{\hbox{$\scriptscriptstyle  R\kern-0.2em  R$}}}}}

\newcommand{\Z}{{{\mathchoice  {\hbox{$\textstyle{\text{\strange Z}}$}}
{\hbox{$\textstyle{\text{\strange Z}}$}}
{\hbox{$\scriptstyle  Z\kern-0.3em  Z$}}
{\hbox{$\scriptscriptstyle  Z\kern-0.2em  Z$}}}}}

\newcommand{\N}{{{\mathchoice  {\hbox{$\textstyle{\text{\strange N}}$}}
{\hbox{$\textstyle{\text{\strange N}}$}}
{\hbox{$\scriptstyle  N\kern-0.3em  N$}}
{\hbox{$\scriptscriptstyle  N\kern-0.2em  N$}}}}}

\renewcommand{\phi}{\varphi}
\usepackage{amsmath,amsthm,amscd}

\begin{document}
\title[Discrete spectrum of quantum layers]{On the discrete spectrum of quantum layers}

 \author{Zhiqin Lu} \address{Department of
Mathematics, University of California,
Irvine, Irvine, CA 92697, USA} \email[Zhiqin Lu]{zlu@uci.edu} 

\author{Julie Rowlett} \address{Max Planck Institut f\"ur Mathematik, Vivatgasse 7,  53111 Bonn, Germany} \email[Julie Rowlett]{rowlett@mpim-bonn.mpg.de}

\begin{abstract}
Consider a quantum particle trapped between a curved layer of constant width built over a complete, non-compact, $\mathcal C^2$ smooth surface embedded in $\mathbb{R}^3$.  We assume that the surface is asymptotically flat in the sense that the second fundamental form vanishes at infinity, and that the surface is not totally geodesic.  This geometric setting is known as {\it a quantum layer.}  We consider the quantum particle to be governed by the Dirichlet Laplacian as Hamiltonian.  Our work concerns the existence of bound states with energy beneath the essential spectrum, which implies the existence of {\it discrete spectrum.}  We first prove that if the Gauss curvature is integrable, and the surface is {\it weakly $\kappa$-parabolic,} then the discrete spectrum is non-empty.  This result implies that if the total Gauss curvature is {\it non-positive,}  then the discrete spectrum is non-empty.  We next prove that if the Gauss curvature is {\it non-negative,} then the discrete spectrum is non-empty.  Finally, we prove that if the surface is {\it parabolic,} then the discrete spectrum is non-empty if the layer is sufficiently thin.  
\end{abstract}  

\maketitle

\pagestyle{myheadings}

\newcommand{\nn}{{\bf n}}
\newcommand{\ka}{K\"ahler }
\newcommand{\ii}{\sqrt{-1}}

\section{Introduction}
Mesoscopic physics describes length scales from one atom to micrometers.  At these scales, the behavior of particles is no longer described by classical physics:  quantum effects are observed.  Numerous phenomena such as quantum dots and wells occur at the scale of mesoscopic physics, and {\it all nanotechnology is on the mesoscopic scale.}   

Consider, for example, electrons trapped between two semi-conducting materials, or more generally, quantum particles trapped between hard walls.  Mathematically, such situations are  described using a {\it quantum layer}. 

Let $p: \Sigma\to\mathbb R^3$ be an embedded surface in $\mathbb R^3$.  We will always make the following assumptions on $\Sigma$.  
\subsubsection{Hypotheses}  \label{hyp} 
\begin{enumerate}
\item $\Sigma$ is a $\mathcal C^2$ smooth surface; 
\item $\Sigma$ is orientable, complete, but non-compact;  
\item $\Sigma$ is not totally geodesic;   
\item $\Sigma$ is asymptotically flat in the sense that the second fundamental form,  denoted $B$ throughout this work, tends to zero at infinity. 
\end{enumerate}

A {\it quantum layer}   over $\Sigma$ is an oriented differentiable manifold $\Omega \cong \Sigma\times [-a,a]$ for some (small) positive number $a$.  Let $\vec N$ be the unit normal vector of $\Sigma$ in $\mathbb R^3$. Define
\[
\tilde p: \Omega\to\mathbb R^3
\]
by 
\[
\tilde p(x,t)=p(x)+t\vec N_x.
\]
If $a$ is small, then $\tilde p$ is clearly an {immersion.} The Riemannian metric $ds_\Omega^2$ is defined as the pull-back of the Euclidean metric via $\tilde p$. The Riemannian manifold $(\Omega, ds^2_\Omega)$ is called the {\it quantum layer}.   Physically, the quantum particles are trapped between two copies of the same semi-conducting material $\Sigma$ at a uniform distance of $2a$ apart, where $a$ is of mesoscopic scale.  A natural question is: 
\begin{quote} 
{ Does there exist a geometric condition on $\Sigma$ which guarantees the existence of bound states {with energy beneath the essential spectrum?}  }
\end{quote}

Let us formulate this precisely.  Let $\Delta=\Delta_\Omega$ be the Laplacian with respect to the Riemannian metric $ds_{\Omega} ^2$, and assume the Dirichlet boundary condition. Since $\Omega$ is a smooth, complete manifold with boundary, the Dirichlet Laplacian is the Friedrichs extension of the Laplacian acting on $\mathcal C_0 ^\infty (\Omega)$ and is self-adjoint.  The spectrum {of the Dirichlet Laplacian} consists of two parts:  discrete, isolated eigenvalues of finite multiplicity and essential spectrum.  We distinguish the eigenvalues which are {\it disjoint} from the essential spectrum and refer to these as the {\it discrete spectrum,} since there may also be embedded eigenvalues within the essential spectrum.   Physically, the quantum particles are governed by the Dirichlet Laplacian as Hamiltonian, { and eigenvalues correspond to the Dirichlet energy of bound states.  Therefore:  
\begin{quote}
{ The existence of discrete spectrum is equivalent to the existence of bound states in the physical model whose energy is beneath the essential spectrum.}  
\end{quote} 
}

Let $\kappa$ be the Gauss curvature of $\Sigma$ throughout this paper.  {Our work is motivated by the following conjecture.}  
\begin{conj}\label{ccon} Under the preceding assumptions \eqref{hyp} on $\Sigma$, if
 \begin{equation}\label{huber}
\int_\Sigma|\kappa|d\Sigma<+\infty,
\end{equation}
then there exists an $\alpha=\alpha(\Sigma)$ such that for all $a \in (0, \alpha)$, the discrete spectrum of the quantum layer over $\Sigma$ of width $2a$ is non-empty.  
\end{conj} 

\begin{rem}  By a theorem of Huber~\cite{huber}, if ~\eqref{huber} is valid, then $\Sigma$ is conformal to a compact Riemann surface with finitely many points removed. Moreover, White~\cite{bwhite}  proved that if 
\begin{equation} \label{white} 
\int_\Sigma \kappa ^-d\Sigma<+\infty, \end{equation} 
where 
$$\kappa = \left\{ \begin{array}{ll} \kappa^+ & \kappa \geq 0 \\ -\kappa^- & \kappa < 0 \end{array} \right., $$
then 
\[
\int_\Sigma |\kappa| d\Sigma<+\infty.
\]
Thus~\eqref{huber} can be weakened to ~\eqref{white}.  
\end{rem}  

Conjecture ~\ref{ccon} was proven under the condition
\begin{equation} \label{dek-neg}
\int_\Sigma \kappa d\Sigma\leq 0
\end{equation} 
through the work of Duclos, Exner and Krej\v{c}i\v{r}{\'\i}k \cite{dek} and Carron, Exner, and Krej\v{c}i\v{r}{\'\i}k \cite{cek-1}.  {Moreover, \cite{cek-1} also proved that Conjecture ~\ref{ccon} holds if the gradient of the mean curvature is locally square integrable, and the total mean curvature is infinite.

Our work focuses on the remaining case:  
\[
\int_\Sigma \kappa d\Sigma>0.
\]
By the theorem of Huber~\cite{huber}, there exists a compact Riemann surface $\bar \Sigma$ and finitely many points $p_1,\cdots, p_s$ such that 
\[
\Sigma=\bar\Sigma\backslash\{p_1,\cdots,p_s\}.
\]
In particular, $\Sigma$ has finitely many ends, $E_i$ ($1\leq i\leq s$).
By a theorem of Hartman~\cite{hartman}*{Theorem 6.1, 7.1}, we have 
\begin{equation}\label{iuy}
\int_\Sigma \kappa d\Sigma  ={2\pi \left( \chi(\Sigma)- \sum_{i=1} ^s \lambda_i \right)}
\end{equation}
where $\chi(\Sigma)$ is the Euler characteristic number of $\Sigma$, and $\lambda_i$ are the isoperimetric constants at each end $E_i$ defined by
$$\lambda_i = \lim_{r\to \infty} \frac{ \textrm{vol} (B(r)\cap E_i)}{\pi r^2}.$$ 
The existence of these limits follows from the integrability of the Gauss curvature \eqref{huber}.  

We also have
$$\chi(\Sigma)\leq \chi(\bar\Sigma)-s=2-2g(\bar\Sigma)-s<2.$$
This together with~\eqref{iuy} implies that $\chi(\Sigma)=1$, and hence the surface is differomorphic to $\mathbb R^2$.  Consequently, 
\begin{equation} \label{int-2p} 0 < \int_\Sigma \kappa d \Sigma \leq 2 \pi. \end{equation} 

Although the topology of the surface is completely known, this is the only remaining case in which the conjecture has not yet been proven.  

We recall the main results of \cites{dek,cek-1}.  
\begin{theorem}[Duclos, Exner and Krej\v{c}i\v{r}{\'\i}k] \label{th:dek}  Let $\Sigma$ be a $\mathcal C^2$-smooth complete simply connected non-compact surface with a pole embedded in $\R^3$.  Let the layer {$\Omega \cong \Sigma \times [-a, a]$} built over the surface be not self-intersecting.  If the surface is not a plane, but it is asymptotically flat, then {if $a$ satisfies condition \eqref{hypa} below, each of the following implies Conjecture ~\ref{ccon}.}    
\begin{enumerate}
\item The Gauss curvature satisfies \eqref{huber} and \eqref{dek-neg};
\item $\Sigma$ is $\mathcal C^3$ smooth, {and $a$ is sufficiently small;}
\item $\Sigma$ is $\mathcal C^3$ smooth, the Gauss curvature is integrable, the gradient of the mean curvature $\nabla_g H $ is $L^2$ integrable, and the total mean curvature is infinite;
\item the Gauss curvature is integrable and $\Sigma$ is cylindrically symmetric.
\end{enumerate} 
\end{theorem}

In \cite{cek-1}, Carron, Exner and Krej\v{c}i\v{r}{\'\i}k proved that the conjecture holds under more general conditions.  {They no longer required the surface to have a pole, and they also proved Conjecture 1 under the additional assumptions that the gradient of the mean curvature is square integrable and the total mean curvature is infinite.}     

\begin{theorem}[Carron, Exner and  Krej\v{c}i\v{r}{\'\i}k]\label{th:cek} Let $\Sigma$ be a complete  asymptotically flat, noncompact connected surface of class $\mathcal C^2$ embedded in $\R^3$ and such that the Gauss curvature satisfies \eqref{huber}.  Let the layer $\Omega$ of width $2a$ be defined so that $\Omega$ does not overlap, and $a$ satisfies condition \eqref{hypa} below.  Then, any of the {following} imply Conjecture ~\ref{ccon}.  
\begin{enumerate}
\item The Gauss curvature satisfies \eqref{dek-neg};
\item $a$ is small enough, and the gradient of the mean curvature $\nabla_g H$ is locally $L^2$ integrable;
\item the gradient of the mean curvature $\nabla_g H$ is $L^2$ integrable, and the total mean curvature is infinite;
\item $\Sigma$ contains a cylindrically symmetric end with a positive total Gauss curvature. 
\end{enumerate}
\end{theorem}   

The general method used in both \cites{dek,cek-1} is:   first compute the infimum of the essential spectrum, next construct appropriate test functions, and finally apply the variational principle to prove that if one of the conditions is satisfied, then there must be an eigenvalue strictly less than the essential spectrum.  The pole and symmetry assumption (4) were necessary in \cite{dek} because their test functions are radially symmetric.  \\

The first main result of the present paper generalizes \cites{dek,cek-1} by demonstrating that Conjecture ~\ref{ccon} holds if the surface is {\it weakly $\kappa$-parabolic.}  (We refer to \S~\ref{4} for the definition of weak $\kappa$-parabolicity.)

\begin{theorem}\label{main6} Let $\Sigma$ be a complete surface in $\R^3$ which satisfies the hypothes{e}s \eqref{hyp}, and assume that $\Sigma$ is weakly $\kappa$-parabolic.  Then, there exists $\alpha > 0$, depending only on the supremum of the norm of the second fundamental form, such that for all $a \in (0, \alpha)$, the {discrete spectrum of the quantum layer over $\Sigma$ of width $2a$ is non-empty.}  
\end{theorem}

Although the proof of Theorem ~\ref{main6} is based on the same principles used in ~\cites{dek, cek-1}, our theorem not only generalizes their results, but also shows that their argument fits nicely into the notion of weak $\kappa$-parabolicity, which provides a geometric abstraction of their argument. \\ 

{The class of layers considered in \cites{cek-1, dek} was already quite broad, but not exhaustive.  For example, a question raised in \cite{dek} which remained unanswered in \cite{cek-1} is: }  
\begin{question} 
{\it {
Does Conjecture ~\ref{ccon} also hold for thick layers built over surfaces of strictly positive total Gauss curvature without assuming cylindrical symmetry or square-integrability of the gradient of the mean curvature? }}
\end{question}  

So-called ``thick layers'' are those  whose thickness satisfies ~\eqref{hypa} and is not otherwise restricted.  The first author and C.~Lin proved in ~\cite{ll-1}*{Theorem 1.1} that Conjecture ~\ref{ccon} holds when $\Sigma$ can be represented as the graph of a convex function satisfying certain conditions.  Our next theorem generalizes \cite{ll-1} {and gives an affirmative answer to the above question under the additional assumption that the Gauss curvature is everywhere non-negative.}   

\begin{theorem}\label{main4}
Let $\Sigma$ be a complete surface in $\R^3$ which satisfies the hypothes{e}s \eqref{hyp}.  Assume that the Gauss curvature of $\Sigma$ is non-negative and satisfies \eqref{huber}.  Then for all $a$ such that
\begin{equation} \label{hypa} a \in (0, B_{\infty} ^{-1}), \quad B_\infty := \sup_{p \in \Sigma} ||B(p)||, \end{equation} 
{the discrete spectrum of} the quantum layer over $\Sigma$ of width $2a$ {is non-empty.}  
\end{theorem}

\begin{rem} The condition ~\eqref{hypa} is merely a technicality to ensure that $\Omega$ is immersed in $\R^3$ and is slightly weaker than the non-overlapping assumption made in \cites{cek-1, dek} (c.f. the Remark on p. 6--7 of \cite{dek}).   
\end{rem} 

The proof of Theorem~\ref{main4} is more subtle.  Test functions similar to those used in \cites{dek,cek-1} rely on the weak $\kappa$-parabolicity of the surfaces, but a surface with non-negative Gauss curvature will {\it not}  be weakly $\kappa$-parabolic.  The main idea is to work on annuli, rather than on disks.  In general, the integration of the mean curvature outside a compact set may be quite small since the surface is asymptotically flat.  But using a result of  White \cite{bwhite}, we actually know that the total mean curvature is at least of linear growth.  This estimate plays a crucial role in the proof.  

In fact, our final main result is more general than Conjecture ~\ref{ccon}.  Any surface whose Gauss curvature is integrable must be {\it parabolic} (see \S 3 for the definition), yet not all parabolic surfaces have integrable Gauss curvature, as demonstrated in \S 5.  

\begin{theorem} \label{parabolic} 
Let $\Sigma$ be a complete, parabolic surface in $\R^3$ which satisfies the hypothes{e}s \eqref{hyp}.  Then, there exists $\alpha > 0$ such that the discrete spectrum of the quantum layer over $\Sigma$ of width $2a$ is non-empty for all $a \in (0, \alpha)$. 
\end{theorem}

This  paper  is organized as follows.  In section \S~\ref{200}, we recall the variational principles for the essential spectrum and the ground state, and we determine the infimum of the essential spectrum.  In \S~\ref{4}, we introduce the notion of $\kappa$-parabolicity and prove Theorem ~\ref{main6}.  The proof of Theorem ~\ref{main4} comprises \S~\ref{400}.   {We conclude in \S~\ref{500} with the proof of Theorem ~\ref{parabolic} and a discussion of further generalizations.}  

\section{Variational principles and the infimum of the essential spectrum}\label{200}
It is well known that
\begin{equation} \label{ground}
\sigma_0=\underset{f\in \mathcal C_0^\infty(\Omega)}{\inf}\frac{\int_\Omega|\nabla f|^2 d\Omega}{\int_\Omega f^2 d\Omega}
\end{equation} 
is the infimum of the spectrum of the  Laplacian.  

{For a compact set $E \subset \Sigma$, we shall use the notation $$f \in \mathcal C^\infty _0 (\Sigma \setminus E)$$ to denote a function $f \in \mathcal C^\infty _0 (\Sigma)$ whose support lies in $\Sigma \setminus E$. The infinimum of the essential spectrum is 
\begin{equation}\label{ess}
\sigma_{ess}=\underset{K}{\sup}\underset{f\in \mathcal C_0^\infty(\Omega\backslash K)}{\inf}\frac{\int_\Omega|\nabla f|^2 d\Omega}{\int_\Omega f^2 d\Omega}
\end{equation}
where $K$ runs over all compact subsets of $\Omega$.  }  
Since $\Omega = \Sigma\times [-a,a]$, it is not hard to see that
\begin{equation}\label{ess-1}
\sigma_{ess}=\underset{K\subset\Sigma}{\sup}\,\,\underset{f\in \mathcal C_0^\infty(\Omega\backslash K\times[-a,a])}{\inf}\frac{\int_\Omega|\nabla f|^2 d\Omega}{\int_\Omega f^2 d\Omega},
\end{equation}
where $K$ runs over all compact subsets of $\Sigma$.

It follows that $\sigma_0 \leq \sigma_{ess}$, and in particular, we have the following.

\begin{rem}\label{ineq} 
If $\sigma_0 < \sigma_{ess}$, then the discrete spectrum is non-empty.  
\end{rem}

Let $(x_1,x_2)$ be a local coordinate system of $\Sigma$. Then $(x_1,x_2,t)$ defines a local coordinate system of $\Omega$.  Such a local coordinate system is called a {\it Fermi coordinate system}. Let $x_3=t$, and let $ds_\Omega^2=G_{ij}dx_idx_j$. Then  we have
\begin{equation}\label{gij}
G_{ij}=\left\{
\begin{array}{ll}
(p+t\vec N)_{x_i}(p+t\vec N)_{x_j} & 1\leq i,j\leq 2;\\
0& i=3, \text {or } j=3, \text{ but } i\neq j;\\
1& i=j=3.
\end{array}
\right.
\end{equation}

We will demonstrate below that the infimum of the essential spectrum is equal to the spectral threshold of the planar quantum layer of width $2a$, namely, $\pi^2/4a^2$.  

We make the following definitions.  For a smooth function $f$ on $\Omega$, let 
\begin{align}
&Q(f,f)=\int_\Omega|\nabla f|^2 d\Omega-\frac{\pi^2}{4a^2}\int_{\Omega} f^2 d\Omega;\label{q}\\
& Q_1(f,f)=\int_\Omega|\nabla' f|^2 d\Omega;\label{q1}\\
& Q_2(f,f)=\int_\Omega\left(\frac{\pa f}{\pa t}\right)^2 d\Omega-\frac{\pi^2}{4a^2}\int_{\Omega} f^2 d\Omega,\label{q2}
\end{align}
where 
$$|\nabla' f|^2=\sum_{i,j=1}^2 G^{ij}\frac{\pa f}{\pa x_i}\frac{\pa f}{\pa x_j}$$
is the square of the  norm of the horizontal differential.

Obviously, we have
\[
Q(f,f)=Q_1(f,f)+Q_2(f,f),
\]
and 
\[
\int_\Omega|\nabla f|^2 d\Omega=\int_\Omega|\nabla' f|^2 d\Omega
+\int_\Omega\left(\frac{\pa f}{\pa t}\right)^2 d\Omega.
\]
 Clearly, we have
\[
\int_\Omega|\nabla f|^2 d\Omega\geq
\int_\Omega\left(\frac{\pa f}{\pa t}\right)^2 d\Omega. 
\]

Let  $ds^2_\Sigma=g_{ij} dx_idx_j$  be the Riemannian metric of $\Sigma$ with respect to the coordinates $(x_1,x_2)$.  We shall compare the matrices $(G_{ij})_{1\leq i,j\leq 2}$ and $(g_{ij})$ outside a big compact set of $\Sigma$.
By~\eqref{gij}, we have
\begin{equation}\label{formula}
G_{ij}=g_{ij}+tp_{x_i}\vec N_{x_j}+tp_{x_j}\vec N_{x_i}+t^2\vec N_{x_i}\vec N_{x_j},
\end{equation}
where we note that $g_{ij} = p_i p_j$.  

Using~\eqref{formula}, we have 
\begin{equation} \label{curv} 
 \det(G_{ij})=\det(g_{ij}) (1-Ht+\kappa t^2).
 \end{equation}
(Note that the mean curvature $H$ is defined to be the trace of the second fundamental form.)   

We assume that at the point $x$, local coordinates $(x_1, x_2)$ are chosen such that $g_{ij}=\delta_{ij}$. We have the estimate 
\begin{equation} \label{G-B}
|G_{ij}-\delta_{ij}|\leq a||B||,
\end{equation} 
where $B$ is the second fundamental form of the surface $\Sigma$.   {Based on these calculations, we have the following.}  

\begin{prop}\label{prop22}
Let $\Sigma$ be an embedded surface in $\R^3$ which satisfies hypotheses \eqref{hyp}.  Then, for any {$a$ which satisfies ~\eqref{hypa},} 
the quantum layer $\Omega \cong \Sigma \times [-a, a]$ is an {immersed} submanifold of $\R^3$.    Moreover,  for any $\eps>0$, there is a compact set $K$ of $\Sigma$ such that on $\Sigma\setminus K$ we have
\begin{equation} \label{g-G} 
(1-\eps)
\begin{pmatrix}
g_{11}&g_{12}\\
g_{21}& g_{22}
\end{pmatrix}
\leq
\begin{pmatrix}
G_{11}&G_{12}\\
G_{21}& G_{22}
\end{pmatrix}
\leq (1+\eps)
\begin{pmatrix}
g_{11}&g_{12}\\
g_{21}& g_{22}
\end{pmatrix}.
\end{equation}
In particular, we have
\begin{equation} \label{dSdO}
(1-\eps)^2d\Sigma dt\leq d\Omega\leq(1+\eps)^2 d\Sigma dt.
\end{equation} 
On the other hand,  there exists $\alpha = \alpha(\Sigma, \eps) > 0$ such that for all $a \in (0, \alpha)$, the above inequalities hold at any point of $\Sigma$. 
\end{prop}
{\bf Proof.}  It follows from ~\eqref{curv} and ~\eqref{G-B} that $\Omega$ is an immersed submanifold whenever $a$ satisfies ~\eqref{hypa}.  By ~\eqref{G-B} and the assumption that $\Sigma$ is asymptotically flat, for any $\epsilon > 0$, there exists a compact subset $K \subset \Sigma$ such that ~\eqref{g-G} holds on $\Sigma \setminus K$.  Since $||B||$ vanishes at infinity, both $H$ and $\kappa$ also vanish at infinity.  Therefore, ~\eqref{dSdO} follows from ~\eqref{curv}.  Finally, for a fixed $\epsilon > 0$, we may choose a compact subset $K \subset \Sigma$ such that both ~\eqref{g-G} and ~\eqref{dSdO} hold on $\Sigma \setminus K$, and since $K$ is compact, by ~\eqref{curv} and ~\eqref{G-B}, we may choose $a$ sufficiently small such that ~\eqref{g-G} and ~\eqref{dSdO} hold also on $K$.  
\qed

The following technical lemma shall be used throughout the remaining sections.  

\begin{lemma} \label{lemmaj}
For a surface $\Sigma$ which satisfies ~\eqref{hyp}, let $ j$ be a smooth function on $\Sigma$ with compact support, and let 
$$\chi(t) = \cos \left( \frac{\pi t}{2a} \right).$$
Then there exist universal constants $C_1$ and $C_2$ which depend only on $B_{\infty}$ (the supremum of the norm of the second fundamental form) and an absolute constant $C_3$  such that 
\begin{equation} 
\label{pqr} 
Q(j \chi t, j \chi t) \leq {a  \int_\Sigma j^2 d\Sigma + C_2 a^3 \int_\Sigma |\n j|^2 d \Sigma+C_3a^3\int_\Sigma j^2\kappa d\Sigma.}
\end{equation}
{Moreover, we also have the following estimate
\begin{equation} \label{rty} Q(j \chi, j \chi) \leq a \int_\Sigma j^2 \kappa d\Sigma + (a+a^2 ||B||_\infty + a^3 ||B||_\infty) \int_\Sigma |\n j|^2 d\Sigma.  \end{equation} }

\end{lemma} 
{\bf Proof.} 
It is a straightforward exercise to compute:   
 \begin{align} \label{calcs}\begin{split} 
 & \int_{-a}^a\chi^2(t) dt=-2\int_{-a} ^a \chi'(t) \chi(t) t dt = a,\\
 & \int_{-a} ^a \chi(t)^2 t^2 dt  =-\frac 23\int_{-a} ^a \chi'(t) \chi(t) t^3 dt = \frac{a^3 (\pi^2 - 6)}{3\pi^2}, \\ 
 &\int_{-a} ^a (\chi'(t))^2 t^2 dt = \frac{a(6+\pi^2)}{12}, \\ 
&\int_{-a} ^a \chi(t)^2 t^4 dt= \frac{a^5 (120 - 20\pi^2 + \pi^4)}{5 \pi^4},\\ 
&\int_{-a} ^a (\chi'(t))^2 t^4 dt = \frac{a^3(20\pi^2-120 +\pi^4)}{20 \pi^2}.
 \end{split}\end{align}
By ~\eqref{curv}, ~\eqref{formula}, and since $j$ is independent of $t$, there is a constant $C_0$ which depends only on $B_{\infty}$ such that 
\begin{align*}
&Q_1 (j \chi t, j \chi t) \leq C_0 \int_{-a} ^a \int_{\Sigma} |\nabla  j|^2 \chi^2 t^2  d\Sigma dt
\end{align*}

By~\eqref{calcs}, we have
\[
Q_1 (j \chi t, j \chi t) 
  \leq C_2a^3\int_\Sigma|\nabla j|^2 d\Sigma, \]
where 
$$C_2 := C_0 \left( \frac{\pi^2 - 6}{3 \pi^2} + \frac{a^2 (120 - 20\pi^2 + \pi^4)}{5 \pi^4} B_\infty\right).$$ 

Next, we have 
\begin{align*} &Q_2 (j \chi t, j \chi t) = \int_{-a} ^a \int_{\Sigma} \left( (\chi')^2 t^2 + \chi^2 + 2 \chi \chi' t \right) (1 {-} Ht + \kappa t^2) j^2 d \Sigma dt\\
&  \qquad  -\frac{\pi^2}{4a^2} \int_{-a}^a\int_\Sigma \chi^2 (1{-} Ht+\kappa t^2) j^2 d\Sigma. 
\end{align*}

Since {$((\chi')^2 t^2 + \chi^2 + 2 \chi \chi' t)$  and $\chi^2 t^2$ are even functions,} we have
\begin{align*} &Q_2 (j \chi t, j \chi t) = \int_{-a} ^a \int_{\Sigma} \left( (\chi')^2 t^2 + \chi^2 + 2 \chi \chi' t \right) (1 + \kappa t^2) j^2 d \Sigma dt\\
&   \qquad  -\frac{\pi^2}{4a^2} \int_{-a}^a\int_\Sigma \chi^2t^2(1+\kappa t^2) j^2 d\Sigma.
\end{align*}

Using~\eqref{calcs}, we have
\begin{equation} \label{lemmaj-est} 
Q_2 (j \chi t, j \chi t) =a  \int_\Sigma j^2 d\Sigma +\left(\frac 43 -\frac{8}{\pi^2}\right)a^3\int_\Sigma j^2\kappa d\Sigma.
\end{equation} 

{To complete the proof of the Lemma, we first compute as above using ~\eqref{calcs} 
$$Q_2(j  \chi, j \chi) = a \int_\Sigma j^2 \kappa d\Sigma.$$
Next, by ~\eqref{G-B}, 
$$Q_1 (j \chi, j \chi) \leq (1 + a ||B||_\infty) \int_\Omega |\n j|^2 \chi^2 d\Omega.$$
By ~\eqref{curv} and since $\chi^2$ is an even function, 
$$Q_1 (j \chi, j \chi) \leq (1+a ||B||_\infty) \int_\Sigma \int_{-a} ^a |\n j|^2 \chi^2 (1 + \kappa t^2) dt d\Sigma.$$
Using ~\eqref{calcs} and since $\kappa \leq ||B||_\infty$, we have
$$Q_1 (j \chi, j\chi) \leq(1+a||B||_\infty) \left( a + \frac{a^3(\pi^2 - 6)}{3 \pi^2} ||B||_\infty \right) \int_\Sigma |\n j |^2 d\Sigma.$$
Finally, estimating the constants and using the hypothesis $a \leq ||B||_\infty ^{-1}$, we have 
$$Q_1 (j \chi, j \chi) \leq \left( a + a^2 ||B||_\infty + a^3 ||B||_\infty\right) \int_\Sigma |\n j|^2 d\Sigma.$$
}

\qed

Based on the preceding results and the variational principle, we are able to determine $\sigma_{ess}$.  The following lemma is originally due to \cites{cek-1, dek}, but we include a short proof for completeness.  

\begin{lemma}\label{lemess}
Let $\Sigma$ be an embedded surface in $\R^3$ which satisfies the hypotheses ~\eqref{hyp}, and assume the Gauss curvature satisfies ~\eqref{huber}.  Then, for any quantum layer $\Omega$ built over $\Sigma$ of width $2a > 0$, where $a$ satisfies ~\eqref{hypa},  
\[
\sigma_{ess} = \frac{\pi^2}{4a^2}.
\]
\end{lemma}

{\bf Proof.} We first prove that 
$$\sigma_{ess} \geq  \frac{\pi^2}{4a^2}.$$
Let $\epsilon >0$ be given, and let $K$ be a compact set of $\Sigma$ as in Proposition ~\ref{prop22}.  {Let $\tilde{K} \subset \Omega$ be defined as 
$$\tilde{K} \cong K \times [-a, a].$$
}
For $f \in C^{\infty} _0(\Omega \setminus {\tilde K})$, by Proposition ~\ref{prop22}, 
\[
\int_{\Omega}\left(\frac{\pa f}{\pa t}\right)^2 d\Omega\geq (1-\eps)^2
\int_\Sigma\int_{-a}^a \left(\frac{\pa f}{\pa t}\right)^2 dtd\Sigma
\geq (1-\eps)^2 \frac{\pi^2}{4a^2}
\int_\Sigma\int_{-a}^a f^2 dtd\Sigma,
\]
where the last inequality follows from the $1$-dimensional Poincar\'e inequality.  By Proposition~\ref{prop22} again, we have 
\[
\int_\Omega |\nabla f|^2 d\Omega\geq\frac{(1-\eps)^2}{(1+\eps)^2}
\frac{\pi^2}{4a^2} \int_\Omega f^2 d\Omega, 
\]
which by the variational principle for $\sigma_{ess}$ implies
\[
\sigma_{ess}\geq \frac{(1-\eps)^2}{(1+\eps)^2}\frac{\pi^2}{4a^2}. 
\]
Letting $\epsilon \to 0$ completes the proof of the first inequality.  

To complete the proof of the lemma, we  demonstrate the estimate 
\begin{equation}\label{lkj}\sigma_{ess} \leq \frac{\pi^2}{4a^2}.\end{equation}
Since the Gauss curvature tends to zero, it is well known (see ~\cite{brooks}) that the infimum of the essential spectrum of $\Sigma$ is zero. Therefore, for any compact set $K$ and any $\eps>0$, there exists a smooth function $\phi\in \mathcal C^\infty_0(\Sigma\setminus K)$ such that
\[
\int_\Sigma|\Delta\phi|^2 d\Sigma\leq\eps\int_\Sigma\phi^2 d\Sigma.
\]
It follows that 
\[
\int_\Sigma|\nabla\phi|^2 d\Sigma=-\int_\Sigma\phi\Delta\phi d\Sigma\leq\sqrt{\int_\Sigma\phi^2 d\Sigma}\cdot \sqrt{\int_\Sigma(\Delta\phi)^2 d\Sigma}\leq\sqrt{\eps}{\int_\Sigma\phi^2 d\Sigma}.
\]

Using Proposition~\ref{prop22}, for sufficiently large $K$, we have
\[
\int_\Omega|\nabla'\phi|^2 d\Omega\leq2(1+\eps)a \int_\Sigma|\nabla\phi|^2 d\Sigma\leq 2\eps(1+\eps)a\int_\Sigma\phi^2 d\Sigma.
\]
As in Lemma ~\ref{lemmaj}, we let 
$$\chi(t)=\cos (\pi t/2a)$$
and consider the function $\phi \chi$ on $\Omega$.     Thus by ~\eqref{calcs} and ~\eqref{curv}, since $\chi^2$ and $(\chi') ^2$ are even functions, we have
 \begin{align*}
&Q_2(\phi \chi,\phi \chi)=\int_\Sigma\phi^2 d\Sigma\int^a_{-a}\left(\left(\chi'(t)\right)^2-\frac{\pi^2}{4a^2}\chi(t)^2\right) dt\\&+\int_\Sigma\phi^2\kappa d\Sigma\int^a_{-a}t^2\left(\left(\chi'(t)\right)^2-\frac{\pi^2}{4a^2}\chi(t)^2\right) dt. \end{align*}
By ~\eqref{calcs},
we have
 \begin{equation}\label{q2phi-chi}
 Q_2(\phi \chi,\phi \chi)=  a\int_\Sigma\phi^2\kappa d \Sigma,
 \end{equation}
 and hence
 \[
 Q(\phi \chi,\phi \chi)\leq 2a\eps(1+\eps)\int_\Sigma\phi^2d\Sigma+ a\int_\Sigma\phi^2\kappa d \Sigma.
 \]

By Proposition~\ref{prop22} and ~\eqref{huber}, \[
\int_\Omega(\phi\chi)^2 d\Omega\geq2(1-\eps_1)a\int_\Sigma\phi^2 d\Sigma,
\]
and
\[
\int_\Sigma\phi^2\kappa d\Sigma\leq\eps_1\int_\Sigma\phi^2 d\Sigma
\]
for some sufficiently small $\eps_1>0$.
 By the variational principle for $\sigma_{ess}$ and the definition of $Q$, we have
 \[
 \sigma_{ess}-\frac{\pi^2}{4a^2}\leq\frac{Q(\phi \chi,\phi \chi)}{\int_\Omega (\phi \chi)^2 d\Omega}\leq 
\frac{2\eps(1+\eps) + \eps_1}{2(1-\eps)}.\]  
This proves the lemma.

 \qed 
 
\section{$\kappa$-parabolicity}\label{4}
 We refer to \cite{li} for the following definition and basic properties of parabolic manifolds. 
 
\begin{definition}
A complete manifold is {\it parabolic}  if it does not admit a positive Green's function.
Otherwise it is {\it non-parabolic}.    
\end{definition}

We first establish the  following well-known result.

\begin{prop}\label{prop123}
Assume that $\kappa\in L^1(\Sigma)$. Then there exists a positive constant $c_1$  such that the volume growth of $\Sigma$ satisfies
\[
\int_{B(R)} d\Sigma<c_1R^2.
\]
\end{prop}

{\bf Proof.}  By the results of Huber~\cite{huber} and Hartman~\cite{hartman}, $\Sigma$ has only finitely many ends, $\{E_i \}_{i=1} ^s$ and at each end,
\[
\lambda_i = \lim_{r\to \infty} \frac{ \textrm{vol} (B(r) \cap E_i)}{\pi r^2}
\]
exists. This proves the proposition.  

\qed

By the above proposition, a surface whose Gauss curvature $\kappa \in L^1$ is a parabolic manifold  (cf.~\cite{li}). 

The definition of parabolicity is equivalent to:  {\it the capacity of any ball of radius $R$ is zero}.   That is, for any positive $R$ and $\epsilon$,  there exists  a smooth function $\phi$ which satisfies
\begin{align*}
&\phi  \in  C^{\infty} _0 {( \Sigma );} \\
& \phi  \equiv  1\text{ on } B(R);\\
&0\leq\phi\leq 1;\\
& \int_\Sigma |\nabla\phi|^2 d\Sigma < \eps. 
 \end{align*} 

For the rest of the paper, we shall repeatedly use the above equivalent capacity definition of  parabolicity.  \\

Parallel to the above, we make the following definition of weak $\kappa$-parabolicity.

\begin{definition} The {\it $\kappa$-capacity} of a subset $E\subset \Sigma$ is defined to be the infimum of 
\[
\int_\Sigma \left( 2 |\nabla\phi|^2+\kappa\phi^2 \right) d\Sigma,
\]
where $\phi$ is a smooth function on $\Sigma$ with compact support, and $\phi\equiv 1$ in a neighborhood of $E$.

We say $\Sigma$ is {\it weakly $\kappa$-parabolic},  if either $\Sigma$ is a minimal {parabolic surface}, or {if} there exists $p\in\Sigma$ such that $H(p)\neq 0$,  and the $\kappa$-capacity of  a neighborhood of $p$ {is non-positive.} \end{definition}

{\bf Proof of Theorem~\ref{main6}.}
By Lemma ~\ref{lemess} and the variational principles, it suffices to prove 
$$\sigma_0<\sigma_{ess}.$$

If $\Sigma$ is a minimal surface, then $\kappa\leq 0$. Consequently the total Gauss curvature  is negative because $\Sigma$ is not totally geodesic.  {In this case using~\eqref{rty} we obtain
\[
Q(\phi\chi, \phi\chi)\leq a\int_\Sigma \kappa\phi^2 d\Sigma+C\int_\Sigma|\nabla\phi|^2 d\Sigma,
\]
where the constant $C$ depends only on $||B||_\infty$.  By the capacity definition of parabolicity, the second term on the right hand side can be made arbitrary small while the first term is negative if the support of $\phi$ is sufficiently large. Thus we conclude that $Q(\phi\chi, \phi\chi)<0$ for a suitable choice of $\phi$, and in this case, the theorem is proven.

To prove the remaining case, let} $p\in\Sigma$ be a point such that $H(p)\neq 0$.  We assume that there is a constant $\eps_1>0$ such that $|H(p)|>\eps_1$ on the ball of radius $\delta$ centered at $p$, which is denoted $B_p(\delta)$, with $\delta$ a fixed positive constant.  For any $\eps_2>0$,  by the  definition of weakly $\kappa$-parabolic (and choosing $\delta > 0$ smaller if necessary), there exists a smooth function $\phi$ with compact support such that $\phi\equiv 1$ on $B_p(\delta)$, and
\begin{equation} \label{eps2}
\int_\Sigma \left( 2|\nabla\phi|^2+\kappa\phi^2 \right) d\Sigma <\eps_2.
\end{equation} 

Let $j$ be a smooth function  such that the support of $j$ is contained in $B_p(\delta)$,  and 
\[
\left|\int_\Sigma Hjd\Sigma\right|>\eps_3>0
\]
for some positive constant $\eps_3>0$. For a suitable choice of orientation, we may assume that
\[\int_\Sigma Hjd\Sigma>\eps_3>0.
\]

By \eqref{rty}, if $a$ is chosen small enough, 
 $$
 Q(\phi\chi,\phi\chi)\leq a\left( \int_\Sigma \left( {2}|\nabla\phi|^2+\kappa\phi^2 \right) d\Sigma\right)\leq a\eps_2.
 $$

By~\eqref{calcs}, we have
 \begin{equation}\label{q7}
 Q(\phi\chi(t),j\chi(t)t)=-\frac{a}{2} \int_\Sigma H j d\Sigma.
 \end{equation}

Let $\epsilon > 0$.   Using all the above estimates {and Lemma ~\ref{lemmaj},} we have 
\[Q(\phi\chi(t)+\eps j\chi(t)t, \phi\chi(t)+\eps j\chi(t)t) < a\eps_2 
- a \eps_3  \eps+c_1 a \eps^2.
\]
Since $\eps_2$  can be made arbitrary small,  and $\eps_3>0$ is fixed, for a suitable choice of $\eps$, we have
$$Q(\phi\chi(t)+\eps j\chi(t)t, \phi\chi(t)+\eps j\chi(t)t) < 0$$
for all $a \in (0, \alpha)$, for $\alpha$ chosen sufficiently small.

\qed

The following sufficient condition for $\kappa$-parabolicity implies that Theorem ~\ref{main6} is indeed a generalization of both  \cites{dek,cek-1}. 

\begin{lemma} Let $\Sigma$ be a complete surface such that
\[
\int_\Sigma |\kappa| d\Sigma<\infty,
\]
and
\begin{equation}\label{k-neg} 
\int_\Sigma \kappa d\Sigma\leq 0.
\end{equation}  
Then $\Sigma$ is weakly $\kappa$-parabolic.
\end{lemma}

{\bf Proof.} First, if the mean curvature $H \equiv 0$ on $\Sigma$, then $\Sigma$ is a minimal surface, and there is nothing to prove.  So, we assume there exists some $p \in \Sigma$ such that $H(p) \neq 0$.  Since $\kappa$ is integrable, $\Sigma$ is parabolic.  That is, for any ball $B(R)$ of radius $R$ centered at $p$, there exists a smooth function $0\leq \phi\leq 1$ with compact support such that $\phi\equiv 1$ on $B(R)$ and 
\[
\int_\Sigma |\nabla\phi|^2 d\Sigma <\eps,
\]
On the other hand, by the integrability of $\kappa$, for sufficiently large $R$,
\[
\int_{\Sigma \backslash B(R)}\kappa\phi^2 d\Sigma<\eps,\quad 
\int_{\Sigma \backslash B(R)}|\kappa| d\Sigma<\eps, 
\]
and therefore, by the assumption~\eqref{k-neg}, 
$$\int_{\Sigma} \kappa \phi^2 d\Sigma < 2\eps.$$
It follows that for any $\eps > 0$ and $R$ sufficiently large, there exists a function $\phi \in \mathcal C^\infty _0 (\Sigma)$ such that $\phi \equiv 1$ on $B(R)$, and 
$$ \int_\Sigma \left( 2|\nabla \phi|^2 + \kappa \phi^2 \right) d \Sigma < 4\eps.$$
{This estimate proves the lemma.}  

\qed  

\begin{rem}
There are many examples of surfaces which are weakly $\kappa$-parabolic but whose total Gauss curvature is positive.  For example, any $\Sigma$ such that 
$$\int_\Sigma \left| \kappa(\Sigma) \right| d\Sigma < \infty, \quad \textrm{and} \quad \int_\Sigma \kappa(\Sigma) d\Sigma < -\eps_0 < 0,$$
is parabolic.  Therefore, choosing $R>0$ sufficiently large, there exists a function $\phi$ such that 
\begin{align*}
&\phi  \in  C^{\infty} _0 {( \Sigma );} \\
& \phi  \equiv  1\text{ on } B(R);\\
&0\leq\phi\leq 1;\\
& \int_\Sigma |\nabla\phi|^2 d\Sigma < \frac{\eps_0}{8},
 \end{align*} 
and 
$$\int_{B(R)} \kappa(\Sigma) d\Sigma < - \frac{1}{2} \eps_0, \quad \int_{B(R)} \kappa(\Sigma) \phi^2 d\Sigma < -\frac{1}{4} \eps_0.$$
In this case, 
$$\int_\Sigma \left( 2\left| \n \phi \right|^2 + \kappa \phi^2 \right) d\Sigma < 0.$$
It follows that the surface is weakly $\kappa$-parabolic.  Now, since $\phi$ is compactly supported, there exists $R' > 0$ such that the support of $\phi$ is contained in $B(R')$.  Then, it is always possible to change $\Sigma$ outside of $B(R')$ such that the volume growth is of order $R$ (for example, we can attach a cylinder $\pa B(R)\times \mathbb R^+$ to the compact manifold $B(R)$) , which by Hartman's result \cite{hartman} implies that for the new $\Sigma'$, 
$$\int_{\Sigma'} \kappa(\Sigma') d\Sigma' = 2\pi > 0.$$
Since 
$$\Sigma \setminus B(R') \cong \Sigma' \setminus B(R'),$$
$\Sigma'$ is still weakly $\kappa$ parabolic, however the total Gauss curvature is positive.
\end{rem}

\section{Proof of Theorem~\ref{main4}}\label{400}
We proceed as in the proof of Theorem ~\ref{main6} by demonstrating that 
$$\sigma_0 < \sigma_{ess}.$$
By Proposition ~\ref{ineq}, this implies that the discrete spectrum is non-empty.  

First, if the Gauss curvature is identically zero, then by ~\cite{ll-4}*{Theorem 2}, the discrete spectrum is non-empty.  
 
Henceforth, we shall assume that there is at least one point of $\Sigma$ at which the Gauss curvature is positive.  Then, by a theorem of Sacksteder~\cite{sack}*{Theorem (*), page 610}, with suitable choice of  orientation, we can assume that the principle curvatures of $\Sigma$ are always nonnegative.  
 
By \eqref{int-2p}, it follows from the results of White \cite{bwhite} ({Theorem 1, p.~ 318}), that there exists an 
$\eps_0>0$
 such that for $R \gg 0$,
 \[
 \int_{\pa B(R)}||B|| >\eps_0,
 \]
 where $B$ is the second fundamental form of $\Sigma$.  Since  the principle curvatures are nonnegative, we have
 \[
 H\geq||B||.
 \]
 Thus we have
 \begin{equation}\label{use}
 \int_{B(R_2)\setminus B(R_1)}H\,d\Sigma\geq  \eps_0 (R_2-R_1)
 \end{equation}
 provided that both $R_1$ and $R_2$ are large enough.

We follow the same general method of \cites{dek, cek-1, ll-1, ll-2, ll-4}.  The main idea is to use the above estimate together with test functions {\it supported in annuli} whose radii tend to infinity. 

Let $\phi\in C_0^\infty(\Sigma \setminus B(\frac R2))$ be a smooth function such that
 \begin{align} \label{phi-90}
 \begin{split}
& {\text{supp}(\phi) \subset B(\frac 52 R) };\\ 
& \phi  \equiv  1  \text{ on } B(2R) \setminus  B(R);  \\ 
&0\leq\phi\leq 1;\\
 &\int_\Sigma|\nabla\phi|^2 d\Sigma  < \eps_1,  
 \end{split}
 \end{align} 
 where  $\eps_1\to 0$ as $R\to\infty$. The existence of such a function $\phi$ is guaranteed by the parabolicity of $\Sigma$.  

Let $\chi$ be defined as in the previous sections.  By Proposition~\ref{prop22}, there is a constant $c_2$ such that 
$$Q_1 (\phi \chi, \phi \chi) \leq c_2a\int_{\Sigma} |\nabla\phi|^2 d\Sigma<c_2a\eps_1.$$ 
Hence
$$Q_1(\phi \chi, \phi \chi) \to 0 \textrm{ as } R \to \infty.$$

Next, we use the same calculations as in the preceding sections to compute 
$$Q_2 (\phi\chi,\phi\chi) = a \int_{\Sigma} \phi^2 \kappa d \Sigma.$$ 
Since the support of $\phi$ is contained in the annulus $B(\frac 52 R) \setminus B(R/2)$, by the integrability assumption \eqref{huber} on $\kappa$, 
$$Q_2 (\phi \chi, \phi \chi) \to 0 \textrm{ as } R \to \infty.$$
Since $Q = Q_1 + Q_2$, there exists $\eps_3 > 0$ such that 
$$Q(\phi\chi,\phi\chi) \leq \eps_3, \quad \eps_3 \to 0 \textrm{ as } R \to \infty.$$
 
Now let's consider a smooth function $j$ on $\Sigma$ with $0\leq j\leq 1$, such that 
\begin{align*}& j \in  C_0^\infty \left(B(\frac 53 R)\setminus B(\frac 43 R)\right); \\
&j \equiv 1 \textrm{ on } B(\frac{19}{12}R) \setminus B(\frac{17}{12}R);  \\ 
&|\nabla j|<2. 
\end{align*}   
We consider the function $j\chi(t)t$.  By Lemma ~\ref{lemmaj}, the integrability of $\kappa$, and since $|\nabla j| <2$, there is an absolute constant $c_1$ such that 
 \[
 Q(j\chi(t)t,j\chi(t)t)\leq c_1 a \int_\Sigma j^2 d\Sigma. 
 \]
 
Next, let's consider $Q(\phi\chi(t),j\chi(t)t)$. Since the support of $j$ is contained in $\{\phi\equiv 1\}$, by~\eqref{q1}, $Q_1(\phi\chi(t),j\chi(t)t)=0$. 
 
The same computation as in~\eqref{q7} shows that 
\[  Q_2(\phi\chi(t),j\chi(t)t)=-\frac{a}{2} \int_\Sigma H j d\Sigma. \] 

  Let $\eps>0$.  By our preceding calculations 
\begin{eqnarray} \label{eq:conj3} 
 Q(\phi\chi(t)+\eps j\chi(t)t, \phi\chi(t)+\eps j\chi(t)t) \nonumber \\ 
 < \eps_3 - \eps a \int_\Sigma H jd\Sigma+\eps^2 c_1a \int_\Sigma j^2 d\Sigma.
\end{eqnarray}

By ~\eqref{use} and the definition of $j$, there is an independent constant $c_1'$ such that 
\[
Q(\phi\chi(t)+\eps j\chi(t)t, \phi\chi(t)+\eps j\chi(t)t)
 <\eps_3 -\frac{a}{6} \eps  R+\eps^2c_1'R^2.
 \]
 
Since $\eps_3 \to 0$ as $R \to \infty$, we may first choose $R$ sufficiently large and then choose $\eps > 0$ appropriately  {so}  that 
for all $a \in (0, B_{\infty} ^{-1})$, 
 \[
 Q(\phi\chi(t)+\eps j\chi(t)t, \phi\chi(t)+\eps j\chi(t)t)
 <0.
 \]

Therefore the discrete spectrum is non-empty.

 \qed

\section{Proof of Theorem ~\ref{parabolic} and Further Discussions}\label{500}
The proof of Theorem ~\ref{main4} can be generalized to demonstrate Theorem ~\ref{parabolic}, which is a {\it stronger} result than Conjecture ~\ref{ccon}.   \\

{\bf Proof of Theorem~\ref{parabolic}.}  We first note that if the mean curvature vanishes identically, 
then since $\Sigma$ is not totally geodesic, $\kappa\not\equiv 0$. Since $H\equiv 0$, {it follows that $\kappa\leq 0$, and therefore the same method in ~\cites{dek, cek-1} can be used to prove the Theorem ~\ref{parabolic}, even without the assumption 
\[
\int_\Sigma (-\kappa) d\Sigma<\infty.
\]
Nonetheless, we include a short proof here for the sake of completeness.   Using the capacity definition of parabolicity, for any $\eps_1 > 0$ there exists a smooth function $\phi$ such that 
\begin{align*}
&\phi  \in  C^{\infty} _0 {( \Sigma );} \\
& \phi  \equiv  1\text{ on } B(2R);\\
&0\leq\phi\leq 1;\\
& \int_\Sigma |\nabla\phi|^2 d\Sigma < \eps_1. 
 \end{align*} 
By~\eqref{rty}, we have
\begin{equation}\label{nmk}
Q(\phi \chi, \phi \chi) \leq aC{\eps_1} + a \int_\Sigma \phi^2 \kappa d \Sigma,\end{equation}
as $R \to \infty$, where $C$ is a constant depending only on $B_\infty$.  Since $\kappa \leq 0$ and $\kappa \not\equiv 0$, for sufficiently large $R$, 
$$Q(\phi \chi, \phi \chi) < 0.$$

Thus, the theorem is proved when the mean curvature is identically zero. 
So from now on we  assume that $H \not \equiv 0$. }  
 
First, we assume the mean curvature $H$ is smooth.  Let $\phi$ be defined as above and let $$j = \phi H.$$
We compute 
$$Q(\phi \chi, \eps j \chi t) \leq \eps c_1 a^2 \int_\Sigma \left| \n \phi \cdot \n j \right| d\Sigma - \eps \frac{a}{2} \int_\Sigma \phi H j d\Sigma, $$
where the constant $c_1$ depends only on $||B||_{\infty}$.  By Lemma \ref{lemmaj} and ~\eqref{pqr} 
$$Q(\eps j \chi t, \eps j \chi t) \leq \eps^2 a \int_\Sigma j^2 d\Sigma + \eps^2 a^3 c_3 \int_\Sigma |\n j|^2 d\Sigma + \eps^2 c_4 a^3 \int_\Sigma j^2 \kappa d\Sigma,$$
where the constant $c_3$ depends only on $||B||_\infty$, and the constant $c_4$ is independent.  

{Therefore, using~\eqref{nmk}, we have 
\begin{align}\label{uyt}
\begin{split}
& 
Q(\phi \chi + \eps j \chi t, \phi \chi + \eps j \chi t) \leq aC\eps_1 \\
&+ a \left( \int_\Sigma \phi^2 \kappa d\Sigma - \eps \int_\Sigma \phi H j d\Sigma + \eps^2 \int_\Sigma j^2 d\Sigma \right)\\& + 2a^2  \eps c_1  \int_\Sigma |\n \phi \cdot \n j | d\Sigma  \\ 
& + a^3 \left( \eps^2 c_3 \int_\Sigma |\n j |^2 d\Sigma + \eps^2 c_4 \int_\Sigma j^2 \kappa d\Sigma \right). 
\end{split}
\end{align} 
}
Note that 
$$\kappa \leq \frac 14H^2.$$
If $\kappa \equiv \frac14 H^2$, then all {points of $\Sigma$ are umbilic. Hence by the Meusnier Theorem (see p. 175 of \cite{m-p}), $\Sigma$ is a sphere, which contradicts the fact that $\Sigma$ is non-compact.   Thus by continuity $\kappa < \frac14 H^2$ on a set of positive measure.  It follows that for} sufficiently large $R$, there exists a fixed $\eps_5 >0$ which depends only on $B$ such that 
{$$ \int_\Sigma \phi^2 \kappa d\Sigma \leq \frac{1}{4} \int_\Sigma \phi^2 H^2 d\Sigma - \eps_5.$$}
Recalling the definition of $j$, we have {
\begin{align}\label{quiz}
\begin{split}
& 
Q(\phi \chi + \eps j \chi t, \phi \chi + \eps j \chi t) \leq aC\eps_1 \\
& + a \left( \left( \frac{1}{4} - \eps+ \eps^2 \right) \int_\Sigma \phi^2 H^2 d\Sigma - \eps_5 \right) + 2a^2 \eps c_1 \int_\Sigma \left| \n \phi \cdot \n (\phi H) \right| d\Sigma \\ 
& + a^3 \left( \eps^2  c_3 \int_\Sigma | \n (\phi H) |^2 d\Sigma + \eps^2 c_4  \int_\Sigma \phi^2 H^2 \kappa d\Sigma \right). 
\end{split}
\end{align} 
}
To make this quantity negative, we first set {
$$\eps = \frac{1}{2}.$$
This implies 
\begin{align}\label{quiz2}
\begin{split}
& 
Q(\phi \chi + j \chi t/2, \phi \chi + j \chi t/2) \leq aC\eps_1  - a  \eps_5  \\ 
& + {a^2 c_1} \int_\Sigma \left| \n \phi \cdot \n (\phi H) \right| d\Sigma \\ 
& + a^3 \left(\frac{c_3}{4} \int_\Sigma | \n (\phi H) |^2 d\Sigma + \frac{ c_4}{4}  \int_\Sigma \phi^2 H^2 \kappa d\Sigma \right). 
\end{split}
\end{align} 

We take $R$ large to make $\eps_1$ small so that $C\eps_1<\frac 12 \eps_5$.  At this point $\phi$, $j$, and $\eps_5$} are fixed, so we can choose $a$ sufficiently small to control the $a^2$ and $a^3$ terms, and therefore, 
$$Q(\phi \chi + j \chi t/2, \phi \chi +  j \chi t/2) < 0.$$

In the general case in which $H$ need not be smooth, there is a canonical smooth approximation,
$$H_a (x) := \frac{1}{\text{ vol}(B_a(x))} \int_{B_a (x)} H d\Sigma,$$
if the underlying space is the Euclidean plane,
where $B_a (x)$ is the geodesic ball of radius $a$ about a point $x \in \Sigma$.  This is an example of what is known in Russian literature as a ``Steklov approximation'' and can be found in \cite{akh} and more recently in \S 4 of \cite{ks}.  On $\Sigma$, we use the mollification technique of \cite{gt}*{Lemma 7.1, 7.2} which, for any $\eps_2 > 0$ guarantees the existence of a smooth function $H_{\eps_1}$ such that 
$$\int_\Sigma |H-H_{\eps_1}|^2 d\Sigma =\eps_3 \to 0 \quad \textrm{ as } \eps_2 \to 0.$$
We replace $j$ defined above by $\phi H_{\eps_1}$ and estimate similarly.   

By~\eqref{uyt}, we have {
\begin{align*} 
\begin{split}
& 
Q(\phi \chi + \eps j \chi t, \phi \chi + \eps j \chi t) \leq aC\eps_1 \\
&+ a \left( \int_\Sigma \phi^2 \kappa d\Sigma - \eps \int_\Sigma \phi H j d\Sigma + \eps^2 \int_\Sigma j^2 d\Sigma \right) \\&+ 2a^2  \eps c_1  \int_\Sigma |\n \phi \cdot \n j | d\Sigma  \\ 
& + a^3 \left( \eps^2 c_3 \int_\Sigma |\n j |^2 d\Sigma + \eps^2 c_4 \int_\Sigma j^2 \kappa d\Sigma \right). 
\end{split}
\end{align*} 
}

{By the discussion above, we assume that  $\kappa < \frac14 H^2$ on a set of positive measure.  Therefore, for sufficiently large $R$, there exists a constant $\eps_5>0$ which depends only on $B$ such that 
$$\int_\Sigma \phi^2 \kappa d\Sigma \leq \frac{1}{4} \int_\Sigma \phi^2 H^2 d\Sigma - \eps_5.$$}
Therefore, 
{\begin{align}\label{quiz3}
\begin{split}
& 
Q(\phi \chi + \eps j \chi t, \phi \chi + \eps j \chi t) \leq aC\eps_1 \\
&+ a \left( \frac{1}{4} \int_\Sigma \phi^2 H^2 d\Sigma - \eps \int_\Sigma \phi H j d\Sigma + \eps^2 \int_\Sigma j^2 d\Sigma - \eps_5 \right) \\ 
&+ 2a^2  \eps c_1  \int_\Sigma |\n \phi \cdot \n j | d\Sigma   + a^3 \left( \eps^2 c_3 \int_\Sigma |\n j |^2 d\Sigma + \eps^2 c_4 \int_\Sigma j^2 \kappa d\Sigma \right). 
\end{split}
\end{align}
}

{Since $j = \phi H_{\eps_1}$, by the Cauchy-Schwarz inequality and the definition of $H_{\eps_1}$, 
$$\left| \int_\Sigma \phi H j d\Sigma - \int_\Sigma \phi^2 H^2 d\Sigma \right| \leq \sqrt{\eps_4} \sqrt{ \int_\Sigma \phi^2 H^2 d\Sigma},$$
where 
$$\eps_4 = \sqrt{ \int_\Sigma \phi^2 (H_{\eps_1} - H)^2 d\Sigma } \to 0 \textrm{ as } \eps_2 \to 0.$$

Moreover, by definition of $H_{\eps_1}$, we may assume that 
$$\int_\Sigma \phi^2 (H_{\eps_1} + H)^2 d\Sigma \leq 9 \int_\Sigma \phi^2 H^2 d\Sigma.$$
It follows from the Cauchy-Schwarz inequality and the definition of $H_{\eps_1}$ that 
$$\left| \int_\Sigma j^2 d\Sigma - \int_\Sigma \phi^2 H^2 d\Sigma \right| \leq 3 \sqrt{\eps_4} \sqrt{\int_\Sigma \phi^2 H^2 d\Sigma}.$$
}

Therefore, by the definition of $j$ we estimate 
\begin{align*}
\begin{split}
& 
Q(\phi \chi + \eps j \chi t, \phi \chi + \eps j \chi t) \leq {aC\eps_1} - a\eps_5 \\
&+ a \left(\left( \frac{1}{4} - \eps + \eps^2 \right) \int_\Sigma \phi^2 H^2 d\Sigma \right) \\
&+ a \eps {\sqrt{\eps_4}\left(1 + 3 \eps \right)} \sqrt{ \int_\Sigma \phi^2 H^2 d\Sigma } \\
&+ 2a^2 \eps c_1  \int_\Sigma |\n \phi \cdot \n j | d\Sigma \\&  + a^3 \left( \eps^2 c_3 \int_\Sigma |\n j |^2 d\Sigma + \eps^2 c_4 \int_\Sigma j^2 \kappa d\Sigma \right). 
\end{split}
\end{align*}

Proceeding similarly, we first let $\eps = \frac{1}{2}$.  Then, we have 
\begin{align*}
\begin{split}
& 
Q(\phi \chi +  j \chi t/2, \phi \chi + j \chi t/2) \leq aC\eps_1 - a\eps_5 \\
&+ \frac{a}{2} {\frac{5 \sqrt{\eps_4}}{2}} \sqrt{ \int_\Sigma \phi^2 H^2 d\Sigma } \\
&+ {a^2 c_1}  \int_\Sigma |\n \phi \cdot \n j | d\Sigma   + a^3 \left( \frac{c_3}{4} \int_\Sigma |\n j |^2 d\Sigma + \frac{c_4}{4} \int_\Sigma j^2 \kappa d\Sigma \right). 
\end{split}
\end{align*}
To make this negative, we first let $R$ be sufficiently large to make 
{$$C \eps_1 < \frac{\eps_5}{2}.$$}   
By the convergence of $H_{\eps_1}$ to $H$, this can be done to make {$\eps_4$} also small.  Finally, since this fixes both $\phi$ and $j$, we can choose $a$ sufficiently small to control the $a^2$ and $a^3$ terms so that 
$$Q(\phi \chi + \eps j \chi t, \phi \chi + \eps j \chi t) < 0.$$ 
Therefore, the discrete spectrum is also non-empty in the non-smooth case, when the mean curvature is merely continuous.  

\qed

\begin{rem} The preceding result implies Conjecture ~\ref{ccon} and is in fact a stronger result, because not all parabolic surfaces have integrable Gauss curvature.  For example, let $f \in \mathcal C^{\infty} (\R)$ satisfy:  
\begin{enumerate}
\item $\int_1^\infty f(r) dr =O(r^2)$;
\item $\int_1^\infty |f''(r)|dr=\infty$.
\end{enumerate}
Let $\Sigma = \R^2$ with polar coordinates $(r, \theta)$, and let the Riemannian metric on $\Sigma$ be $g = dr^2+f(r)^2 d\theta^2$.  Then  the volume growth is quadratic, and hence $(\Sigma, g)$ is parabolic (see \cite{li}), but by (2) above, the Gauss curvature is not integrable.  
\end{rem}

\subsection{Further Discussions}  

In \S 7 of  \cite{dek}, they construct a layer whose discrete spectrum is empty, but this example is not asymptotically planar.  In both our work and \cites{cek-1, dek}, in some cases we are able to prove the existence of discrete spectrum for ``thick layers,'' meaning those whose thickness satisfies ~\eqref{hypa} and is not further restricted:  Theorem ~\ref{main4}; Theorem ~\ref{th:dek} (1), (3), (4); Theorem ~\ref{th:cek} (1), (3), (4).  However, for the remaining cases, we must assume the layer is sufficiently thin to prove the existence of discrete spectrum.  It would be interesting to investigate whether Theorem ~\ref{parabolic} holds for all ``thick layers'' over parabolic surfaces $\Sigma$ satisfying ~\eqref{hyp}, or whether one may construct a layer of width $2a$ for $a$ satisfying ~\eqref{hypa} over a parabolic surface which satisfies ~\eqref{hyp}, such that the discrete spectrum 
of the layer is empty.  The results of \cites{cek-1, dek} and Theorem ~\ref{main6} indicate that we expect the following conjecture holds.    

\begin{conj} \label{con2} Let $\Sigma$ be a surface which satisfies ~\eqref{hyp}, and assume the Gauss curvature satisfies ~\eqref{huber}.  Then, there exists $\alpha = \alpha(B_\infty) > 0$ depending only on the supremum of the norm of the second fundamental form such that for all $a \in (0, \alpha)$, the discrete spectrum of the quantum layer of width $2a$ over $\Sigma$ is non-empty.  
\end{conj} 

Although we are unable to completely prove the above conjecture, the following results imply the conjecture holds under certain geometric assumptions. 

\begin{theorem}\label{main5}
Let $\Sigma$ be a complete surface in $\R^3$ which satisfies the hypotheses ~\eqref{hyp}, and assume that the Gauss curvature of $\Sigma$ satisfies ~\eqref{huber}.  
The following is sufficient to imply Conjecture ~\ref{con2}:  
 \begin{eqnarray} \label{assu} 
\exists \eps_0 > 0 \textrm{ such that } \forall \, R>0, \, \exists j \in \mathcal  C^\infty _0 (\Sigma \setminus B(R)) \textrm{ with } \\ 
\left(\int_\Sigma Hjd\Sigma\right)^2>\eps_0\int_\Sigma (|\nabla j|^2+j^2) d\Sigma. \nonumber
\end{eqnarray}
\end{theorem}

{\bf Proof.}   As noted previously, we may assume 
$$\int_\Sigma \kappa d\Sigma > 0.$$
Since $\Sigma$ is parabolic, for any $R>0$  sufficiently large there exists a smooth function $\phi\in C_0^\infty(B(R')\setminus B(R))$ such that
 \begin{align} 
 \begin{split}
& \phi  \equiv  1  \text{ on } {B(R'/2)} \setminus  B(2R);  \\ 
 &\int_\Sigma|\nabla\phi|^2 d\Sigma  < \eps_1;\\
 &0\leq\phi\leq 1;\\
 &\int_\Sigma \kappa\phi^2 d\Sigma<\eps_1,
 \end{split}
 \end{align} 
 where $R'$ is a sufficiently large number.
It follows from the arguments in the proof of Theorem~\ref{main4} that 
\[
Q(\phi\chi,\phi\chi)<\eps_3, \quad\eps_3\to 0 \text{ as } R\to\infty.
\]
We choose $R'$ large enough such that 
\[
{\rm supp}\,(j)\subset B(R'/2)\setminus B(2R).
\]
Without loss of generality, we assume that
\[
\int_\Sigma Hjd\Sigma>0.
\]
Then using the same method as in the proof of Theorem~\ref{main4}, letting $\eps>0$, we have
\begin{align*}
&
Q(\phi\chi+\eps j\chi(t) t, \phi\chi+\eps j\chi(t) t)\\
&\leq \eps_3- a \eps \int_\Sigma Hjd\Sigma+c_1\eps^2 a \int_\Sigma (j^2 + |\nabla j|^2)d\Sigma
\end{align*}
for some constant $c_1$ which is {\it independent of $R$}.    Since $\eps_3\to 0$ as $R\to\infty$,  we first choose $R$ sufficiently large to make $\eps_3$ sufficiently small.  Then, using \eqref{assu}, for all sufficiently small $\eps > 0$, 
$$Q(\phi\chi+\eps j\chi(t) t, \phi\chi+\eps j\chi(t) t) < 0.$$
This implies that the discrete spectrum of the layer $\Omega = \Sigma \times [-a,a]$ is non-empty, for any $a$ which satisfies ~\eqref{hypa}.

\qed 

Based on the above theorem, we make the following {\it purely Riemannian geometric}  conjecture.   

\begin{conj} \label{con3}Let $\Sigma$ be a complete surface in $\R^3$   which satisfies the hypotheses ~\eqref{hyp}, and assume that the Gauss curvature $\kappa$ of $\Sigma$ satisfies ~\eqref{huber}.  If the total Gauss curvature is positive 
$$\int_\Sigma \kappa > 0,$$ 
then ~\eqref{assu} holds.\end{conj}  

\begin{rem}  By the results of \cites{cek-1, dek}, Conjecture ~\ref{con3} together with Theorem ~\ref{main5} would imply Conjecture~\ref{con2}. \end{rem}

\begin{rem} It is straightforward to verify that conditions (3) and (4) in Theorem ~\ref{th:cek} are each implied by ~\eqref{assu}.  For example, to prove that (3) implies \eqref{assu}, $j$ is replaced by $jH$, and ~\eqref{assu} follows from a direct calculation.  Assuming (4), ~\eqref{assu} follws from either \cite{dek}*{Lemma 6.1} or \cite{cek-1}*{page 783}. \end{rem} 

The following proposition uses the result of White \cite{bwhite} to demonstrate a weaker  version of the inequality of ~\eqref{assu}.  

\begin{prop} Assume that $\Sigma$ is a complete surface in $\R^3$ which satisfies the hypotheses ~\eqref{hyp}, and that the Gauss curvature $\kappa$ satisfies
$$\int_\Sigma  |\kappa|  d \Sigma < \infty, \quad \int_\Sigma \kappa d \Sigma > 0.$$
Then for any $R>0$ there exists a positive constant $\eps_0$ and a function $j\in C_0^\infty( \Sigma \setminus B(R) )$ which satisfies 
\begin{equation}\label{assu-2}
\left(\int_\Sigma j |H| d\Sigma\right)^2 > \eps_0\int_\Sigma (|\nabla j|^2+j^2) d\Sigma.
\end{equation}
\end{prop}

{\bf Proof.}   By the result ~\cite{bwhite}*{Theorem 1, p.~318}, for sufficiently large $R$,
\[
\int_{\pa B(R)} ||B||>c_1>0
\] 
for some constant $c_1$.  

By the parabolicity of $\Sigma$ and Proposition~\ref{prop123}, we can find a function $j$ whose support is contained in  $B(\frac 53R) \setminus B(\frac 43 R)$, with 
\begin{align}
\begin{split} 
&j \equiv  1 \text{ on } B(\frac{19}{12}R) \setminus B(\frac{17}{12}R); \\
&\quad\int_\Sigma |\nabla j|^2 d\Sigma <1,\\
&\quad\int_\Sigma j^2 d\Sigma \leq  c_1 R^2,
\end{split}\end{align} 
where $c_1$ is a  positive constant independent of $R$.  Hence we have
\[
\int_\Sigma ||B||jd\Sigma\geq c_2 R
\]
for some constant $c_2$. 

On the other hand, since $\kappa$ is integrable, for sufficiently large $R$, we have 
\[
\int_\Sigma j\sqrt{|\kappa|} d\Sigma\leq\sqrt{\int_\Sigma jd\Sigma}\cdot\sqrt{\int_\Sigma j\kappa d\Sigma}\leq \eps_3{R},
\]
for some small positive constant $\eps_3$ when $R$ is large. 

Since $|H|\geq||B||-\sqrt{2|\kappa|}$, the above inequalities show that 
$$\int_{\Sigma} j |H| d \Sigma \geq \int_\Sigma ||B||j d\Sigma - \int_\Sigma j \sqrt{2 |\kappa|} d\Sigma \geq (c_2  - \sqrt 2\eps_3){R}.$$
Therefore, there exists $\eps > 0$ such that for sufficiently large $R$,
$$\int_{\Sigma} j |H| d \Sigma > \eps R.$$

Finally, by definition of $j$, for sufficiently large $R$, 
$$\int_\Sigma( |\nabla j|^2 + j^2) d\Sigma < 2 c_1 R^2.$$

Therefore, there exists a constant $\eps_0 > 0$ such that for all $R$ sufficiently large, 
$$\left(\int_{\Sigma} j |H| d \Sigma \right)^2 > \eps_0  \int_\Sigma (|\nabla j|^2+j^2) d\Sigma.$$

\qed 

\begin{rem}
Note that the estimate ~\eqref{assu-2} is weaker than ~\eqref{assu}, because the integration is $j |H|$ rather than $jH$.  Although we are unable with our present methods to prove Conjecture ~\ref{con3}, the proposition supports the conjecture since it shows that  if the mean curvature has fixed sign off some compact set, then Conjecture ~\ref{con3} holds.  However, by our methods, we cannot prove the conjecture when the mean curvature $H$ continually oscillates between positive and negative all the way to infinity; c.f. Example 6 from \S 6 of \cite{dek}. One would need a new and different argument to prove Conjecture ~\ref{con3} in such cases. 
\end{rem}

Our final theorem below shows that if the surface $\Sigma$ satisfies certain isoperimetric inequalities, this is sufficient for ~\eqref{assu}.    

\begin{theorem}
Let $\Sigma$ satisfy the hypotheses \eqref{hyp}, and assume $\Sigma$ also satisfies the following.   
\begin{enumerate}
\item
The isoperimetric inequality holds. That is, there is a positive constant $\delta_1$ such that if $D$ is a domain  in $\Sigma$, we have
\[
(length(\pa D))^2\geq \delta_1 Area(D).
\]
\item There is another positive constant $\delta_2$ such that for any compact set $K \subset \Sigma$, there is a curve $C \subset \Sigma \setminus K$ such that if $\vec \gamma$ is the  normal vector of $C$  in $\Sigma$, then there is a vector $\vec \eta$ in $\mathbb R^3$ such that
\[
\langle \vec \gamma,\vec \eta \rangle\geq\delta_2>0.
\]
\item All such curves $C$ are tamed. That is, let $\sigma(t,x)$ be the geodesic flow of $\vec\gamma$. Then there exist constants $\delta_3, A$ such that the following hold.
\begin{enumerate}
\item[(3.1)] $\sigma(t,x)$ is defined up to $|t|<\delta_3$; 
\item[(3.2)] the map $C\times (-\delta_3,\delta_3)\to \Sigma$ is diffeomorphic onto its image;  
\item[(3.3)] the derivatives of $\sigma$ and its inverse are bounded by the fixed constant $A$. 
\end{enumerate}
\end{enumerate} 
Then, if the Gauss curvature satisfies \eqref{huber}, ~\eqref{assu} is valid.  
\end{theorem}

{\bf Proof.}   
To construct the required function $j$ satisfying~\eqref{assu}, we let $\rho$ be a smooth non-increasing cut-off function such that 
$$
 \rho(t) = \left\{ \begin{array}{ll} 1, & t \leq 1; \\ 0, & t \geq 1. \end{array}  \right . 
$$
Let $C$ be the curve outside a compact set $K$, and let $D$ be the compact domain of $\Sigma$ such that $\pa D=C$. Assume that $\vec\gamma$ is the outward norm of $C$ in $\Sigma$.

Define the cut-off function $\tilde\rho$ on $\Sigma$ as follows
\[
\tilde\rho=\left\{
\begin{array}{ll}
1& x\in D\\
\rho(\eps^{-1}dist(x,C))& x\not\in D
\end{array}
\right.,
\]
where $\eps$ is a positive number to be determined later.

Let $(x,y,z)$ be the standard coordinates of $\mathbb R^3$.
 Without loss of generality, assume the vector $\vec \eta$ in the hypotheses of the theorem is the $z$-direction in three-dimensional Euclidean space. 

 Let $\vec n$ be the normal vector of $\Sigma\subset \mathbb R^3$. 
Let $n_z$ be the $z$-component of $\vec n$.  Define the function $j$ by
\[
j=\tilde\rho\cdot n_z.
\]

Since  $|\rho n_z|+|\nabla(\rho n_z)|\leq ||B||_\infty+1$, 
if we choose $\eps$ small enough, we have
\[
\int_\Sigma(|\nabla j|^2+j^2)d\Sigma\leq C\eps^{-2} Area(D)
\]
for $\eps$ small, where the constant $C$ depends on {$A$.}  

On the other hand, since  $Hn_z=\Delta z$, we have
$$
 \left|\int_\Sigma H\tilde\rho n_z d\Sigma\right| =\left|\int_\Sigma\nabla z\nabla \tilde\rho d\Sigma \right| $$ $$ =(\langle\vec\gamma,\nabla z\rangle+o(1))\cdot length(C)>\frac 12\delta_2 \cdot length(C). $$
Therefore using the isoperimetric inequality, the conclusion of the theorem holds.
 
 \qed

\section*{Acknowledgements} 
Both authors are grateful to the anonymous referee whose comments improved the quality of this paper.    The first author is partially supported by  NSF grant DMS-09-04653.  The second author gratefully acknowledges the support of the Max  Planck Institut f\"ur Mathematik in Bonn, Germany.  Finally, both authors would like to thank D.~Krej{\v{c}}i{\v{r}}{\'{\i}}k for helpful comments which enriched the content of this paper, and in particular his suggestion for the test function used in the proof of Theorem ~\ref{parabolic} and for drawing our attention to references \cites{ks, akh}.

\end{document}